\documentclass[10pt]{article}

\usepackage{amsmath} 
\usepackage[mathscr]{eucal}
\usepackage{mathrsfs}
\usepackage{amssymb} 
\usepackage{theorem}
\usepackage{rotating}
\usepackage{stmaryrd}

\usepackage{hyperref}
\usepackage[capitalise]{cleveref}

\theoremstyle{change} 
\newtheorem{theorem}{Theorem}[section] 
\newtheorem{lemma}[theorem]{Lemma} 
\newtheorem{proposition}[theorem]{Proposition}
\newtheorem{corollary}[theorem]{Corollary}
\theorembodyfont{\rmfamily} 
\newtheorem{remark}[theorem]{Remark}
\newtheorem{example}[theorem]{Example}

\newtheorem{notation}[theorem]{Notation}
\newtheorem{nothing}[theorem]{} 

\newenvironment{proof}{\noindent{\bf Proof}\ }{\qed\bigskip}

\renewcommand{\le}{\leqslant}
\renewcommand{\ge}{\geqslant} 
\renewcommand{\leq}{\leqslant}

\newcommand{\Aut}{\mathrm{Aut}}

\newcommand{\betatilde}{\tilde{\beta}}
\newcommand{\BIGOP}[1]
  {\mathop{\mathchoice
  {\raise-0.22em\hbox{\huge $#1$}}
  {\raise-0.05em\hbox{\Large $#1$}}{\hbox{\large $#1$}}{#1}}}
\newcommand{\bl}{\mathrm{bl}}
\newcommand{\Bl}{\mathrm{Bl}}

\newcommand{\catfont}{\mathsf}

\newcommand{\dia}{\mathrm{diag}}

\newcommand{\FF}{\mathbb{F}}

\newcommand{\GL}{\mathrm{GL}}
\newcommand{\Hom}{\mathrm{Hom}}

\newcommand{\id}{\mathrm{id}}
\newcommand{\idem}{\mathrm{idem}}
\newcommand{\Ind}{\mathrm{Ind}}

\newcommand{\lexp}[2]{\setbox0=\hbox{$#2$} \setbox1=\vbox to
                 \ht0{}\,\box1^{#1}\!#2}
\newcommand{\lmod}[1]{\llap{\phantom{|}}_{#1}\catfont{mod}}

\newcommand{\Mat}{\mathrm{Mat}}

\newcommand{\qed}{\nobreak\hfill
                   \vbox{\hrule\hbox{\vrule\hbox to 5pt
                   {\vbox to 8pt{\vfil}\hfil}\vrule}\hrule}}
\newcommand{\QQ}{\mathbb{Q}}

\newcommand{\Res}{\mathrm{Res}}

\newcommand{\rk}{\mathrm{rk}}

\newcommand{\ZZ}{\mathbb{Z}}

\title{The ring of perfect $p$-permutation bimodules for blocks with cyclic defect groups\footnote{{\bf MR Subject Classification:}  
20C20, 19A22{\bf Keywords:}  $p$-permutation modules; trivial source modules; blocks of group algebras; idempotents.}}

\author{\small Robert Boltje\\
  \small Department of Mathematics\\
  \small University of California\\
  \small Santa Cruz, CA 95064\\
  \small U.S.A.\\
  \small boltje@ucsc.edu
  \and
  \small Nariel Monteiro\thanks{Research was supported by the NSF grant No. 2213166}\\
  \small Department of Mathematics\\ 
  \small University of California\\
  \small Santa Cruz, CA 95064\\
  \small U.S.A.\\
  \small namontei@ucsc.edu}
\date{August 7, 2024}

\begin{document}
\sloppy


\maketitle




\begin{abstract}
Let $B$ be a block algebra of a group algebra $FG$ of a finite group $G$ over a field $F$ of characteristic $p>0$. This paper studies ring theoretic properties of the representation ring $T^\Delta(B,B)$ of perfect $p$-permutation $(B,B)$-bimodules and properties of the $k$-algebra $k\otimes_\ZZ T^\Delta(B,B)$, for a field $k$. We show that if the Cartan matrix of $B$ has $1$ as an elementary divisor then $[B]$ is not primitive in $T^\Delta(B,B)$. If $B$ has cyclic defect groups we determine a primitive decomposition of $[B]$ in $T^\Delta(B,B)$. Moreover, if $k$ is a field of characteristic different from $p$ and $B$ has cyclic defect groups of order $p^n$ we describe $k\otimes_\ZZ T^\Delta(B,B)$ explicitly as a direct product of a matrix algebra and $n$ group algebras.
\end{abstract}


\section{Introduction}\label{sec intro}

Throughout this paper $F$ denotes a field of characteristic $p>0$ and $G$ denotes a finite group. We assume that $F$ contains a root of unity of order $|G|_{p'}$, or equivalently, that $F$ is a splitting field for the group algebras $FH$ of all subgroups $H\le G$.

\smallskip
The trivial source ring $T(FG)$, or ring of $p$-permutation $FG$-modules, plays a particular role in the modular representation theory of finite group (see Section~\ref{sec T(FG)} for precise definitions of $T(FG)$ and and related representation rings). In particular, if $H$ is another finite group, $p$-permutation $(FG,FH)$-bimodules that are projective on each side feature prominently in notions of equivalences between block algebras, as for instance splendid Rickard equivalence (see \cite{rickard1996splendid}), splendid Morita equivalence (see \cite[Section 9.7]{Li18b}), $p$-permutation equivalence (see \cite{boltjexu2008p-permutation}, \cite{linckelmann2009trivial}, \cite{boltje2020ppermutation}), functorial equivalence (see \cite{boucyilmaz2022diagonal}), and related conjectures, as for instance Brou\'e's abelian defect group conjecture (see \cite{broue1990isometries}) and Puig's finiteness conjecture (see \cite[Conjecture~6.2]{Li18b}).
Its representation group is denoted by $T^\Delta(FG,FH)$. If $G=H$, the group $T^\Delta(FG,FG)$ has a ring structure induced by taking tensor products of bimodules. The classes $[B_i]$ of the block algebras $B_1,\ldots, B_t$ of $FG$, viewed as $(FG,FG)$-bimodules, form a set of pairwise orthogonal idempotents in $T^\Delta(FG,FG)$ whose sum is equal to the identity $[FG]$ of $T^\Delta(FG,FG)$. In a conversation, more than a decade ago, Jacques Th\'evenaz raised the question if the (in general non-central) idempotents $[B_i]$ are primitive in $T^\Delta(FG,FG)$, and if not, how they decompose further. This is equivalent to asking if and how the identity element $[B]$ of the ring $T^\Delta(B,B)$ decomposes further into pairwise orthogonal idempotents.

\smallskip
In this paper we give partial answers to this questions. First note that $T^\Delta(B,B)$ contains the representation group $Pr(B,B)$ of projective $(B,B)$-bimodules as an ideal. The following theorem shows that in general $[B]$ is not primitive in $T^\Delta(B,B)$. More precisely:

\bigskip\noindent
{\bf Theorem A}\quad {\it Let $B$ be a block algebra of $FG$ and let $r$ be the multiplicity of the elementary divisor $1$ of the Cartan matrix of $B$. Then there exists a set of $r$ pairwise orthogonal idempotents in $Pr(B,B)$. Moreover, $r$ is maximal with this property.
}

\bigskip\noindent
Besides the ring $T^\Delta(B,B)$, also its finite-dimensional $k$-algebra version $k T^\Delta(B,B)$, for a field $k$, is of interest. If $k$ has characteristic $0$, it appears as the endomorphism algebra of a functor attached to $B$ in \cite{boucyilmaz2022diagonal}. Its ideal $k Pr(B,B)$ has a particularly nice structure when the characteristic of $k$ is different from $p$.

\bigskip\noindent
{\bf Theorem B}\quad {\it Let $B$ be a block algebra of $FG$ and $k$ a field of characteristic different from $p$. Then the ideal $k Pr(B,B)$ of the $k$-algebra $k T^\Delta(B,B)$ has an identity element which is a central idempotent of $k T^\Delta(B,B)$. Moreover, $k Pr(B,B)$ is isomorphic to $\Mat_l(k)$ as $k$-algebra, where $l$ is the number of isomorphism classes of simple $B$-modules. In particular, $k Pr(B,B)$ is a block algebra of $k T^\Delta(B,B)$.
}

\bigskip\noindent
If $B$ has cyclic defect groups, we obtain complete and explicit (see also Example~\ref{ex special C}) answers about idempotent decompositions in $T^\Delta(B,B)$ and the algebra structure of $kT^\Delta(B,B)$ when the characteristic of $k$ is different from $p$ (see the proof of Theorem~D in Section~\ref{sec proof of Thm C} for explicit formulas).

\bigskip\noindent
{\bf Theorem C}\quad {\it Let $B$ be a block algebra of $FG$ with cyclic defect groups and let $l$ be the number of isomorphism classes of simple $B$-modules. Then there exist pairwise orthogonal primitive idempotents $e_1,\ldots,e_l$ of $T^\Delta(B,B)$ whose sum equals $1_{T^\Delta(B,B)}=[B]$. Moreover, the ring $T^\Delta(B,B)/Pr(B,B)$ is commutative, its only idempotents are $0$ and $1$, and if $f_1,\ldots, f_r$ are pairwise orthogonal idempotents of $T^\Delta(B,B)$ then $r\le l$ and at most one of the elements $f_1,\ldots, f_r$ does not lie in $Pr(B,B)$. Finally $0$ and $1$ are the only central idempotents in $T^\Delta(B,B)$.
}

\bigskip\noindent
{\bf Theorem D}\quad {\it Let $k$ be a field of characteristic different from $p$ and let $B$ be a block algebra of $FG$ with a cyclic defect group $D$ of order $p^n$ and inertial quotient $E$. Then, for $i\in\{0,\ldots,n\}$, the ideal $k\otimes_\ZZ T^\Delta_{\le i}(B,B)$ of $k T^\Delta(B,B):=k\otimes_\ZZ T^\Delta(B,B)$, spanned over $k$ by classes of indecomposable modules with vertex of order at most $p^i$, has an identity element $e_i$. Moreover, the elements $e_0, e_1-e_0, e_2-e_1, \ldots, e_n-e_{n-1}$ form a set of pairwise orthogonal idempotents of $Z(k\otimes_\ZZ T^\Delta(B,B))$ whose sum is the identity element $e_n$ of $k T^\Delta(B,B)$. The $k$-algebra $k T^\Delta(B,B) e_0$ is isomorphism to $\Mat_{|E|}(k)$ and, for $i\in \{1,\ldots,n\}$, the $k$-algebra $T^\Delta(B,B)(e_i-e_{i-1})$ is isomorphic to the group algebra $k[\Aut(D_i)/\pi_i(E)\times \Hom(E,F^\times)]$, where $D_i$ denotes the subgroup of $D$ of order $p^i$ and $\pi_i(E)$ is the subgroup of the abelian group $\Aut(D_i)$ of conjugation automorphisms with elements from $E$. In particular, one has an isomorphism
\begin{equation}\label{eqn decomposition}
   kT^\Delta(B,B) \cong \Mat_{|E|}(k) \times \prod_{i=1}^n k[\Aut(D_i)/\pi_i(E)\times \Hom(E,F^\times)]
\end{equation}
of $k$-algebras. 

Further, for any field $k$, independent of its characteristic, the $k$-algebra $kT^\Delta(B,B)$ is semisimple if and only if $|\Aut(D)|$ is invertible in $k$.
}

\bigskip\noindent
If $B$ has a trivial defect group then $T^\Delta(B,B)=Pr(B,B)\cong \ZZ$ as rings and Theorems  A--D hold for trivial reasons.

\medskip\noindent
The paper is arranged as follows. In Section~\ref{sec T(FG)} we introduce the trivial source ring $T(FG)$ and related rings associated to bimodules, in particular the ring $T^\Delta(FG,FG)$. In Section~\ref{sec idem of Pr(B,B)} we study idempotents in $Pr(B,B)$ and prove Theorems~A and B, see Theorem~\ref{thm Pr(B,B)}(b), Theorem~\ref{thm QPr(B,B)} and Corollary~\ref{cor QPr(B,B)}. 
The rest of the paper is dedicated to proving Theorem~C. 
In Section~\ref{sec normal Sylow} we explicitly describe all indecomposable $p$-permutation $FG$-modules in the case that $G$ has a normal Sylow $p$-subgroup, see Theorem~\ref{thm normal Sylow} which is of independent interest. Section~\ref{sec bouc} recalls the construction of an extended tensor product of bimodules and a theorem that describes tensor product of two induced bimodules via a Mackey formula, both due to Bouc. In the cyclic defect group case, the ring $T^\Delta(B,B)$ is isomorphic to the ring $T^\Delta(F[D\rtimes E], F[D\rtimes E])$, where $D$ is a defect group of $B$ and $E$ the inertial quotient, see Theorem~\ref{thm reduction} which uses the existence of a splendid Rickard equivalence between $B$ and its Brauer correspondent due to Rouquier, and the structure theory of blocks with normal defect group due to K\"ulshammer. Therefore, we explicitly determine the ring structure of $T^\Delta(F[D\rtimes E], F[D\rtimes E])$ in Section~\ref{sec special G}, see Proposition~\ref{prop mult formulas}. Section~\ref{sec gen idem} gives two general results on idempotents that are used in Section~\ref{sec proof of Thm C} for the proofs of Theorem~C and D.

\bigskip\noindent
{\bf Notation and Convention}\quad
If $A$ and $B$ are sets, $A\subseteq B$ indicates that $A$ is a subset of $B$ and $A\subset B$ indicates that $A$ is a proper subset of $B$.

{\em Rings} are not necessarily unital. The set of idempotents of a ring $R$ is denoted by $\idem(R)=\{e\in R\mid e^2=e\}$. Thus, $0$ is an idempotent of $R$. Two idempotents $e,e'$ of $R$ are called {\em orthogonal} if $ee'=e'e=0$. An idempotent $e$ of $R$ is called {\em primitive} in $R$ if $e$ cannot be written as the sum of two non-zero orthogonal idempotents of $R$. A {\em decomposition} of an idempotent $e$ of $R$ is a finite set $\{e_1,\ldots, e_r\}$ of pairwise orthogonal idempotents of $R$ whose sum is equal to $e$. If additionally each $e_i$, $i=1,\ldots,r$, is primitive, we call $\{e_1,\ldots,e_r\}$ e primitive decomposition of $e$. The empty set is considered as a {\em primitive decomposition} of $0$. For any ring $R$, we denote the center of $R$ by $Z(R)$. If $R$ is unital, its unit group is denoted by $R^\times$. 

For any positive integer $l$, any $i,j\in\{1,\ldots,l\}$ and any unital ring $R$, we denote by $E_{ij}\in\Mat_l(R)$ the $l\times l$-matrix with $1$ as the $(i,j)$-entry and zero everywhere else. 

All modules are left modules unless otherwise stated. For any ring $R$, we denote the category of finitely generated left $R$-modules by $\lmod{R}$.

For a subgroup $H$ of $G$ and $g\in G$ we set $\lexp{g}{H}:=gHg^{-1}$. For any a commutative ring $k$, we denote by $kG$ or $k[G]$ the  group algebra of $G$ over $k$. If $M$ is a left $kH$-module and $g\in G$ then $\lexp{g}{M}$ denotes the left $k[\lexp{g}{H}]$-module with underlying $k$-module $M$, but with the action $ghg^{-1} \cdot m= hm$ for $h\in H$ and $m\in \lexp{g}{M}$.

\bigskip\noindent
{\bf Acknowledgement}\quad Parts of this paper were established while we were visiting the University of Valencia, City University London, and the University of Manchester. We would like to express our gratitude for the hospitality experienced in these Mathematics Departments. The second author was supported by the National Science Foundation MPS-Ascend Postdoctoral Research Fellowship under Grant No. 2213166.


\section{The trivial source ring and its variants}\label{sec T(FG)}

We fix a finite group $G$ and a field $F$ of characteristic $p>0$ such that $F$ is a splitting field of $FH$ for all subgroups $H\le G$.
Recall that a {\em $p$-permutation $FG$-module} is a finitely generated $FG$-module $M$ such that for every $p$-subgroup $P$ of $G$ there exists a $P$-stable $F$-basis of $M$, i.e., the restriction of $M$ to $P$ is a {\em permutation $FP$-module}. Equivalently, $M$ is isomorphic to a direct summand of a permutation $FG$-module. Further equivalently, every indecomposable direct summand of $M$ has the trivial module as source. Indecomposable $p$-permutation modules are also called {\em trivial source modules}. Projective modules and permutation modules are examples of $p$-permutation modules.

\smallskip
For any idempotent $e$ of $Z(FG)$, we denote the free abelian group on the set of isomorphism classes $[M]$ of indecomposable $p$-permutation $FGe$-modules $M$ by $T(FGe)$. If $M$ is a (not necessarily indecomposable) $p$-permutation $FGe$-module and $M\cong M_1\oplus\cdots\oplus M_n$ with indecomposable $FG$-modules $M_1,\ldots,M_n$, then we set $[M]:=[M_1]+\cdots+[M_n]\in T(FGe)$. By the Krull-Schmidt Theorem, the modules $M_1,\ldots, M_n$ are uniquely determined up to order and isomorphism. Thus, the classes $[M]$ of indecomposable $p$-permutation $FG$-modules form a finite $\ZZ$-basis of $T(FGe)$. The abelian group $T(FG)$ (for $e=1$) is a unital commutative ring with multiplication $[M]\cdot [M']:=[M\otimes_F M']$. The identity element is the class $[F]$ of the trivial $FG$-module $F$. The span of the classes of indecomposable projective $FGe$-modules in $T(FGe)$ is denoted by $Pr(FGe)$. Note that $Pr(FG)$ is an ideal of $T(FG)$.

\smallskip
Let also $H$ be a finite group. We always identify the $F$-algebras $F[G\times H]$ and $FG\otimes_F FH$ via $(g,h)\mapsto g\otimes h$. If $e \in Z(FG)$ and $f\in Z(FH)$ are idempotents then $e\otimes f$ is an idempotent in $Z(FG\otimes_F FH) = Z(F[G\times H])$. 
If $M$ is an $(FG,FH)$-bimodule then $M$ can be regarded as left $F[G\times H]$-module via $(g,h)m=gmh^{-1}$, and vice-versa. Under this identification, $(FGe,FHf)$-bimodules can be viewed as $F(G\times H)(e\otimes f^*)$-modules, where $-^*\colon FH\to FH$ is the $F$-linear extension of $h\mapsto h^{-1}$. We set $T(FGe,FGf):=T(F[G\times H](e\otimes f^*))$. The subgroup $T^\Delta(FGe,FHf)$ is defined as the $\ZZ$ span of the classes of indecomposable $p$-permutation $F[G\times H](e\otimes f^*))$-modules which have {\em twisted diagonal} vertices. These are subgroups of $G\times H$ of the form $\Delta(P,\phi,Q):=\{(\phi(y),y)\mid y\in Q\}$, where $P\le G$ and $Q\le H$ are $p$-subgroups and $\phi\colon Q\to P$ is an isomorphism. Equivalently $M$ is projective on either side when regarded as $(FGe,FHf)$-bimodule. Bimodules which are projective on both sides are also called {\em perfect}. Further, we denote by $Pr(FGe,FHf)$ the subgroup of $T^\Delta(FGe,FHf)$ spanned over $\ZZ$ by the classes of projective indecomposable $F[G\times H](e\otimes f^*)$-modules. 

\smallskip
If also $K$ is a finite group then the tensor product $-\otimes_{FH}-$ induces a $\ZZ$-bilinear map $-\cdot_H -\colon T(FG,FH)\times T(FH,FK)\to T(FG,FK)$ such that $T^\Delta(FG,FH)\cdot_H T^\Delta(FH,FK)\subseteq T^\Delta(FG,FK)$, see for instance Theorem~\ref{thm bouc}. 
In particular, if $G=H=K$ and $e\in Z(FG)$ is an idempotent then $T(FGe,FGe)$ is a unital ring with identity element $[FGe]$. 
We always view $T(FGe,FGe)$ as a subring of $T(FG,FG)$. But note that these two rings have different identity elements. By the above, $T^\Delta(FGe,FGe)$ is closed under multiplication and therefore forms a unital subring of $T(FGe,FGe)$ with the same identity element, since $[FGe]\in T^\Delta(FGe,FGe)$. 
In fact, $FGe$ is a direct summand of the permutation $(FG,FG)$-bimodule $FG$ and $FG$ is $\Delta(G)$-projective as $F[G\times G]$-module, where $\Delta(G):=\{(g,g)\mid g\in G\}$, since, after identifying $G$ and $\Delta(G)$ via $g\mapsto (g,g)$, one has an isomorphism $FG \cong \Ind_{\Delta(G)}^{G\times G}(F)$ of $F[G\times G]$-modules.
Moreover, $Pr(FGe,FGe)$ is contained in $T^\Delta(FGe,FGe)$ and is an ideal of $T^\Delta(FGe,FGe)$.

\smallskip
Let $n$ be maximal such that there exists a perfect indecomposable $p$-permutation $F[G\times G](e\otimes e^*)$-module with a vertex of order $p^n$. For $i\in\{0,\ldots,n\}$, let $T_i^\Delta(FGe,FGe)$ (resp.~$T_{\le i}^\Delta(FGe,FGe)$) denote the $\ZZ$-span of the standard basis elements $[M]$ of $T^\Delta(FGe,FGe)$ where $M$ has a vertex of order $p^i$ (resp.~at most $p^i$). Then, 
\begin{equation*}
   T^\Delta(FGe,FGe)= \bigoplus_{i=0}^n T_i^\Delta(FGe,FGe)
\end{equation*}
and, by Theorem~\ref{thm bouc}, $T_{\le i}^\Delta(FGe,FGe)$ is an ideal of $T^\Delta(FGe,FGe)$. Note that $T^\Delta_0(FGe,FGe)=Pr(FGe,FGe)$.

\smallskip
A {\em block idempotent of $FG$} is a primitive idempotent of $Z(FG)$. We denote by $\bl(FG)$ the set of block idempotents $b$ of $FG$ and by $\Bl(FG)$ the set of their corresponding {\em block algebras} $B=FGb$. Thus, $\bl(FG)$ is a primitive decomposition of $1$ in $Z(FG)$, $FG=\bigoplus_{B\in\Bl(FG)} B$ is a decomposition of $FG$ into indecomposable $F[G\times G]$-modules, and $\{[B]\mid B\in\Bl(FG)\}$ is a decomposition of $1=[FG]$ in $T^\Delta(FG,FG)$. In this paper we study if and how the idempotent $[B]$ of $T^\Delta(FG,FG)$ decomposes further into pairwise orthogonal idempotents in $T^\Delta(FG,FG)$, or equivalently, in $[B] \cdot_G T^\Delta(FG,FG) \cdot_G [B]=T^\Delta(B,B)$.


\section{Idempotents in $Pr(B,B)$}\label{sec idem of Pr(B,B)}

Let $F$ and $G$ be as in Section~\ref{sec T(FG)}. Moreover, let $e\in Z(FG)$ be an idempotent and set $B:=FGe$. In this section we will describe maximal sets of pairwise orthogonal primitive idempotents in the ideal $Pr(B,B)$ of $T^\Delta(B,B)$. As special case for $e=1_{FG}$ (resp.~$e\in\bl(FG)$) we obtain $B=FG$ (resp. $B\in\Bl(FG))$. 

\begin{notation}\label{not Rc}
We start with an elementary construction for arbitrary rings. Given a ring $R$ and an element $c$ of $R$, we define the ring $R_c$ to have the same elements and addition as $R$ and multiplication defined by
$$a\cdot_c b=acb\,,$$ 
where the multiplication on the right is the multiplication in $R$. Even if $R$ is unital, $R_c$ in general is not.
The following lemma is a straightforward verification, left to the reader.
\end{notation}

\begin{lemma}\label{lem Rc=Rd}
Let $c$ and $d$ be elements of a unital ring $R$ such that $c=udv$, where $u$ and $v$ are units of $R$. Then, the map $f\colon R_c \to R_d$, $r\mapsto vru$, is a ring isomorphism ring.
\end{lemma}

\begin{nothing}\label{noth Pij}
For the rest of this section we fix the following situation. Choose pairwise orthogonal primitive idempotents $\epsilon_1,\ldots, \epsilon_l$ of $B$. Then $P_i:=FG\epsilon_i$, $i=1,\ldots,l$, is a complete list of pairwise non-isomorphic projective indecomposable $B$-modules. 
Thus, $l=l(B)$ is the number of isomorphism classes of simple $B$-modules.
Every projective indecomposable $F[G\times G](e\otimes e^*)$-module, when viewed as $(FGe,FGe)$-bimodule, is of the form $P_{ij}:=FG\epsilon_i\otimes_F \epsilon_jFG$ with $i,j\in\{1,\ldots,l\}$. 
Moreover, they are pairwise non-isomorphic. Thus, their classes form a $\ZZ$-basis $[P_{ij}]$ of $Pr(FGe,FGe)$. 
Let $C=(c_{ij})\in\Mat_l(\ZZ)$ denote the {\em Cartan matrix} of $B$, i.e., $c_{ij}=\dim_F \epsilon_iFG\epsilon_j= \dim_F \epsilon_iB\epsilon_j$. 
\end{nothing}

\begin{lemma}\label{lem right to left}
Let $\epsilon$ be an idempotent of $FG$. Then, the right $FG$-module $\epsilon FG$, when viewed as left $FG$-module by $g\cdot a:= ag^{-1}$ for $a\in \epsilon FG$ and $g\in G$, is isomorphic to the left $FG$-module $FG\epsilon^*$. In particular, with the above notation, one has $P_{ij}\cong FG\epsilon_i\otimes_F FG\epsilon_j^*$, when $P_{ij}$ is viewed as left $F[G\times G]$-module.
\end{lemma}

\begin{proof}
The map $-^*\colon FG\to FG$, $g\mapsto g^{-1}$, restricts to an isomorphism $\epsilon FG\to FG\epsilon^*$ of left $FG$-modules, where the domain is equipped with the left $FG$-module structure from the statement.
\end{proof}

\begin{lemma}\label{lem Cartan}
{\rm (a)} With the above notation one has $P_{ij}\otimes_{FG} P_{rs}= P_{is}^{c_{jr}}$, for all $i,j,r,s\in\{1.\ldots,l\}$.

\smallskip
{\rm (b)} The map $[P_{ij}]\mapsto E_{ij}$, for $i,j\in\{1,\ldots,l\}$, defines a ring isomorphism between $Pr(B,B)$ and $\Mat_l(\ZZ)_C$, where $C$ is the Cartan matrix of $B$.
\end{lemma}

\begin{proof}
(a) This follows from the following isomorphisms of $(FG,FG)$-bimodules:
\begin{equation*}
   P_{ij} \otimes_{FG} P_{rs} \cong FG\epsilon_i\otimes \bigl(\epsilon_jFG\otimes_{FG} FG \epsilon_r)\bigr) \otimes \epsilon_s FG \cong  P_{is}^{\dim_F \epsilon_jFG\epsilon_r}\,.
\end{equation*}

\smallskip
(b) This follows immediately from (a).
\end{proof}

Recall that the Smith normal form of the Cartan matrix $C$ is a diagonal matrix $D=\dia(d_1,\ldots,d_l)$ with $p$-powers $d_1,\ldots, d_l$ satisfying $d_1|d_2 \cdots | d_l$, see Lemma~6.31 in \cite[Chapter~3]{nagao2014representations}. The numbers $d_1,\ldots,d_l$ are the invariant factors and at the same time the elementary divisors of $C$.

\begin{theorem}\label{thm Pr(B,B)}
{\rm (a)} There exists a ring isomorphism $Pr(B,B)\cong \Mat_l(\ZZ)_D$, where $D\in\Mat_l(\ZZ)$ is the diagonal matrix whose diagonal entries $d_1,\ldots,d_l$ are the elementary divisor of $C$ in ascending size.

\smallskip
{\rm (b)} Let $r$ be the multiplicity of 1 as elementary divisor of the Cartan matrix of $B$. Then there exists a set of $r$ pairwise orthogonal primitive idempotents of $Pr(B,B)$. Moreover, any set of pairwise orthogonal primitive idempotents in $Pr(B,B)$ has cardinality at most $r$.
\end{theorem}

\begin{proof}
Since there exist $U,V\in\GL_l(\ZZ)$ with $D=UCV$, Part~(a) follows from Lemma~\ref{lem Rc=Rd} and \ref{lem Cartan}.

By Part~(a), it suffices to show the statement in (b) for the ring $\Mat_l(\ZZ)_D$ instead of $Pr(B,B)$. First note that, by elementary matrix computations, the $r$ elements $\{E_{11},\ldots,E_{rr}\}$ form a set of $r$ pairwise orthogonal idempotents of $\Mat_l(\ZZ)_D$. They are also primitive, since $E_{ii}\cdot_D \Mat_l(\ZZ)_D \cdot_D E_{ii}= \ZZ d_i^2 E_{ii}\cong \ZZ$ is indecomposable as abelian group. 

Next suppose that $E_1,\ldots, E_s$ are pairwise orthogonal idempotents of $\Mat_l(\ZZ)_D$, set $E=E_1+\ldots+E_s$, and consider the subgroup $M:=E\cdot \ZZ^n$ of $\ZZ^n$. Since $EDE=E$, we obtain $ED\cdot M = M$. This implies that, also after reduction modulo $p$, the map $\bar{E}\bar{D} \colon \bar{M}\to\bar{M}$, $\bar{m}\mapsto \bar{E}\bar{D}\bar{m}$, is surjective. Thus, $\rk_\ZZ M = \dim_{\FF_p}\bar{M}\le \rk_{\FF_p} \bar{D} = r$. On the other hand, we claim that $s\le \rk_\ZZ M$ which then implies $s\le r$ as desired. For $i=1,\ldots, s$, consider the non-zero subgroup $M_i:= E_i\cdot \ZZ^n=(E\cdot_D E_i)\cdot\ZZ^n = E\cdot ((DE_i)\cdot \ZZ^n)\subseteq M$. It suffices to show that the sum of these subgroups of $M$ is a direct sum. Let $x_1,\ldots, x_s\in\ZZ^n$ with $E_1\cdot x_1+\cdots +E_s\cdot x_s=0$. Then, for any $i=1\ldots, s$, multiplication with the matrix $E_iD$ yields $0=E_iDE_i\cdot x_i=E_i\cdot x_i$ and we are done.
\end{proof}

\begin{remark}\label{rem computing idempotents}
One can explicitly determine a maximal set of pairwise orthogonal primitive idempotents of $Pr(B,B)$ by following the proof of Theorem~\ref{thm Pr(B,B)}. In fact, with the notation of the proof, one can apply the inverses of the isomorphism $Pr(B,B)\cong \Mat_l(\ZZ)_C\cong \Mat_l(\ZZ)_D$ to the elements $E_{11},\ldots,E_{rr} \in \Mat_l(\ZZ)_D$.
\end{remark}

The situation for the ring $k Pr(B,B):=k\otimes_{\ZZ}Pr(B,B)$ and over base fields of characteristic different from $p$ is much simpler.

\begin{theorem}\label{thm QPr(B,B)}
Suppose that $k$ is a field of characteristic different from $p$.
There exists a multiplicative $k$-linear isomorphism $k Pr(B,B)\to \Mat_l(k)$. In particular, $k Pr(B,B)$ is a $k$-algebra isomorphic to $\Mat_l(k)$, any primitive decomposition of $1_{k Pr(B,B)}$ has cardinality $l=l(B)$, and any two primitive idempotents of $k Pr(B,B)$ are conjugate by a unit of $k Pr(B,B)$. More explicitly, one has $1_{k Pr(B,B)}=\sum_{i,j=1}^n c'_{ij} [P_{ij}]$, where $(c'_{ij})=C^{-1}\in\Mat_l(k)$, and the elements $\epsilon_i:=\sum_{j=1}^l  c'_{ij} [P_{ij}]$, $i=1,\ldots,l$, form a primitive decomposition of $1_{\QQ Pr(B,B)}$. \end{theorem}

\begin{proof}
By Lemma~\ref{lem Cartan}, mapping the $k$-basis elements $[P_{ij}]$ to $E_{ij}$ induces a $k$-linear isomorphism $k Pr(B,B)\to \Mat_l(k)_C$. Since $\det(C)$ is a power of $p$, $C$ is invertible in $\Mat_l(k)$. Thus, Lemma~\ref{lem Rc=Rd} with $u:=C$ and $v:=I_l$, the identity matrix, yields an isomorphism $\Mat_l(k)_C\to \Mat_l(k)_{I_l}=\Mat_l(k)$. Following the explicit isomorphisms establishes the formula for the identity element of $k Pr(B,B)$ and the primitive decomposition of the identity element corresponding to the primitive idempotents $E_{ii}$, $i=1,\ldots,l$, of $\Mat_l(k)$. 
The remaining statements follow from well-known properties of $\Mat_l(k)$.
\end{proof}

\begin{corollary}\label{cor QPr(B,B)}
If $k$ is a field of characteristic different from $p$
then the identity element of $k Pr(B,B)$ is a central idempotent of $k T^\Delta(B,B):=k\otimes T^\Delta(B,B)$.
Moreover, $Pr(B,B)$ has an identity element if and only every block algebra direct summand of $B$ has defect $0$.
\end{corollary}

\begin{proof}
We may assume that $k$ has at least three elements by passing to an extension field. Then, by \cite[Theorem~1.13.7]{Li18a}, the units in $k T^\Delta(B,B)$ generate $k T^\Delta(B,B)$ as a $k$-vector space. Thus, we only need to show that $1_{k Pr(B,B)}$ is fixed under conjugation by any unit of $k T^\Delta(B,B)$. Since $k Pr(B,B)$ is an ideal in $k T^\Delta(B,B)$, conjugation by a unit induces by a ring automorphism of $k Pr(B,B)$, and therefore must preserve its identity element.

$Pr(B,B)$ has an identity element if and only if $1_{\QQ \Pr(B,B)}\in Pr(B,B)$. Since the elements $[P_{ij}]$ form a $\ZZ$-basis of $Pr(B,B)$ and by the explicit formula for $1_{\QQ Pr(B,B)}$ in Theorem~\ref{thm QPr(B,B)}, one has $1_{\QQ Pr(B,B)}\in Pr(B,B)$ if and only if $\det(C)\in\{-1,1\}$. This happens if and only if all the elementary divisors of $C$ are equal to one. Since the Cartan matrix of a block with defect group $D$ has at least one elementary divisor $|D|$, see Theorem~6.35 in \cite[Chapter~3]{nagao2014representations}, the latter is equivalent to $B$ being a direct sum of block algebras of $FG$ of defect $0$.
\end{proof}

\begin{example}\label{ex special C}
In this example we assume that $b$ is a block idempotent with Cartan matrix 
\begin{equation}\label{eqn special C}
C = I_l+mE_l=  
\begin{pmatrix}
   m +1 & m & \cdots & m\\
   m & m+1 & \ddots &\vdots \\
   \vdots & \ddots & \ddots &m \\
   m & \cdots & m &m+1 
\end{pmatrix} \,, 
\end{equation}
where $m\ge 1$, $I_l\in\Mat_l(\ZZ)$ is the identity matrix, and $E_l\in\Mat_l(\ZZ)$ is the matrix all of whose entries are equal to $1$. This happens for instance when $B=FGb$ has a cyclic and normal defect group, see \cite[Theorem~17.2]{alperin1986local}. In this case $B$ is a Brauer tree algebra whose Brauer tree is a star with exceptional vertex of multiplicity $m$ in the center. We will examine this case further in Section~\ref{sec special G}.

We compute the Smith normal form of $C$. First subtracting the $l$-th row from all the other rows, then subtracting $m$ times the $i$-th row from the $l$-th row, for $i=1,\ldots, l-1$, and finally adding the the $i$-th column to the $l$-th column for $i=1,\ldots,l-1$, results in the diagonal matrix $\dia(1,\ldots,1,ml+1)$, which therefore is the Smith normal form of $C$. 
Furthermore, if $D$ denotes a defect group of $B$, by Theorem~6.35 in \cite[Chapter~3]{nagao2014representations}, the largest elementary divisor of $C$ is $|D|$. This implies that $|D|=ml+1$.

Following Remark~\ref{rem computing idempotents} by computing the matrices $U,V\in\GL_l(\ZZ)$ with $D=UCV$ which result from the above elementary row and column operations, we obtain that the elements $\epsilon_i:=[P_{ii}]-[P_{li}]$, $i=1,\ldots,l-1$, form a maximal set of pairwise orthogonal primitive idempotents of $Pr(B,B)$. This can also be checked directly using Lemma~\ref{lem Cartan}, once one has found these elements.

Now let $k$ be a field of characteristic different from $p$. One quickly verifies that $C^{-1} = I_l-\frac{m}{|D|}E_l\in\Mat_l(k)$. Thus, by Theorem~\ref{thm QPr(B,B)}, one has
\begin{equation*}
   1_{k Pr(B,B)}= \sum_{i=1}^l [P_{ii}] - \frac{m}{|D|}\sum_{i,j=1}^l [P_{ij}]\,.
\end{equation*}
and obtains a primitive decomposition $\{\epsilon_1,\ldots,\epsilon_l\}$ of $1_{k Pr(B,B)}$ by computing 
\begin{equation*}
   \epsilon_l = 1_{k Pr(B,B)}-(\epsilon_1 +\cdots+\epsilon_{l-1})= \sum_{i=1}^l [P_{li}] -\frac{m}{|D|}\sum_{i,j=1}^l [P_{ij}]\,.
\end{equation*}
\end{example}


\section{Trivial source modules for groups with normal Sylow $p$-subgroup}\label{sec normal Sylow}

\begin{theorem}\label{thm normal Sylow}
Let $G$ be a finite group with a normal Sylow $p$-subgroup $P$, let $Q\le P$, set $H:=N_G(Q)$, and let $C$ be a complement of the normal Sylow $p$-subgroup $N_P(Q)$ of $H$. 
 Further, let $S$ be a simple $FC$-module and view it as a simple $F(QC)$-module via the canonical isomorphism $C\cong QC/Q$. Then $M_S:=\Ind_{QC}^G(S)$ is an indecomposable $p$-permutation $FG$-module with vertex $Q$. Moreover, the map $S\mapsto M_S$ induces a bijection between the set of isomorphism classes of simple $FC$-modules and the set of isomorphism classes of indecomposable $p$-permutation $FG$-modules with vertex $Q$.
\end{theorem}

\begin{proof}
The theorem will follow from the following two claims together with the Green correspondence.

{\em Claim  1:}\quad If $S$ is a simple $FC$-module then $U_S:=\Ind_{QC}^H(S)$ is an indecomposable $p$-permutation $FH$-module with vertex $Q$, and $S\mapsto U_S$ induces a bijection between the set of isomorphism classes of simple $FC$-modules and the set of isomorphism classes of indecomposable $p$-permutation $FH$-modules with vertex $Q$.

{\em Claim 2:}\quad If $S$ is a simple $FC$-module then $\Ind_H^G(U_S)$ is indecomposable.

\smallskip
To prove Claim~1, note that we have canonical isomorphisms $C\cong QC/Q \cong H/N_P(Q)$. Therefore, the simple $FH$-modules and the simple $F(QC)$-modules can be identified with the simple $FC$-modules, and it suffices to show that the module $\Ind_{QC/Q}^{H/Q} (S)$ is a projective cover of  the simple $F(H/Q)$-module $S$. To see this, let $S$ and $T$ be simple $FC$-modules and consider
\begin{equation*}
   \Hom_{F(H/Q)}(\Ind_{QC/Q}^{H/Q} (S), T) \cong \Hom_{F(QC/Q)}(S,\Res^{H/Q}_{QC/Q}(T)) 
   =\Hom_{F(QC/Q)}(S,T)\,.
\end{equation*}
If $S\cong T$ then the latter is isomorphic to $F$ and if $S\not\cong T$ then it is isomorphic to $\{0\}$. This shows that the radical quotient of $\Ind_{QC/Q}^{H/Q}(S)$ is isomorphic to $S$. Moreover, $\Ind_{QC/Q}^{H/Q}(S)$ is projective as $QC/Q$-module, since $S$ is projective as $QC/Q$-module. This proves Claim~1.

To prove Claim 2, restrict $\Ind_H^G(U_S)=\Ind_{QC}^G(S)$ to $P$ and apply Mackey's decomposition formula to obtain
\begin{equation}\label{eqn Mackey 1}
   \Res^G_P(M_S)=\bigoplus_{g\in [G/PC]} \Ind_{P\cap\lexp{g}{Q}}^P \Res^{\lexp{g}{Q}}_{P\cap\lexp{g}{Q}}(\lexp{g}{S}) 
   \cong \bigoplus_{g\in [G/PC]} \Ind_{\lexp{g}{Q}}^P (F^{\dim_F(S)})\,.
\end{equation}
The latter expression is a direct sum of indecomposable $FP$-modules with vertices of order $|Q|$. On the other hand, by the Green correspondence, $M_S$ decomposes into one indecomposable direct summand with vertex $Q$ and possibly others with vertex of order smaller than $|Q|$. But such other summands can't exist, since every indecomposable direct summand of $\Res^G_P(M_S)$ has vertices of order $|Q|$. This proves Claim~2 and the proof of the theorem is complete.
\end{proof}


\section{Direct product groups and tensor products of bimodules}\label{sec bouc}

In this section, we recall a construction and a theorem due to Bouc, see \cite{bouc2010inducedbimodules}. Throughout this section, $G$, $H$ and $K$ denote finite groups and $R$ denotes a commutative ring.

\medskip
For a subgroup $X\le G\times H$, we set
\begin{equation*}
   k_1(X):=\{g\in G\mid (g,1)\in X\}\quad\text{and}\quad k_2(X):=\{h\in H\mid (1,h)\in X\}\,.
\end{equation*}
These are subgroups of $G$ and $H$, respectively, and normal subgroups of the images $p_1(X)$ and $p_2(X)$ under the projection maps $p_1\colon G\times H\to G$ and $p_2\colon G\times H\to H$, respectively.
For subgroups $X\le G\times H$ and $Y\le H\times K$ we set
\begin{equation*}
   X*Y:=\{(g,k)\in G\times K \mid \exists h\in H\colon (g,h)\in X\text{ and } (h,k)\in Y\}\,.
\end{equation*}
This is a subgroup of $G\times K$. An element $h\in H$ as above is called a {\em connecting element} for $(g,k)\in X*Y$. We also set $k(X,Y):=k_2(X)\cap k_1(Y)\le H$. If $M$ is a left $RX$-module and $N$ is a left $RY$-module then $M$ can be considered as a right $k(X,Y)$-module via $mh:=(1,h^{-1})m$ for $m\in M$ and $h\in k(X,Y)$ and $N$ can be considered as left $k(X,Y)$-module via $hn:=(h,1)n$ for $h\in k(X,Y)$ and $n\in N$. The resulting $R$-module $M\otimes_{Rk(X,Y)} N$ has the structure of an $R[X*Y]$-module via
\begin{equation*}
   (g,k)(m\otimes n):= (g,h)m\otimes (h,k)n\,,
\end{equation*}
for $(g,k)\in X*Y$, where $h\in H$ is a connecting element for $(g,k)\in X*Y$. The expression on the right hand side of the above definition does not depend on the choice of $h$, since $h$ is unique up to left multiplication with elements in $k(X,Y)$. Just for the purpose of this paper, we denote this $R[X*Y]$ module by $M*N$, suppressing the dependence on $X$ and $Y$. For more properties of the functor $-*-\colon \lmod{RX}\times \lmod{RY}\to\lmod{R[X*Y]}$, see \cite[Section~6]{boltje2020ppermutation}.

\smallskip
In the sequel we will use the following useful properties of the construction $(X,Y)\mapsto X*Y$, whose elementary proof is left to the reader. Recall that for subgroups $P\le G$ and $Q\le H$ and for an isomorphism $\alpha\colon Q\to P$ we set $\Delta(P,\alpha,Q):=\{(\alpha(x),x)\mid x\in Q\}$. We call groups of this type {\em twisted diagonal subgroups} of $G\times H$.

\begin{lemma}\label{lem omnibus}
{\rm (a)} For any $X\le G\times H$, $Y\le H\times K$, $g\in G$, $h\in H$ and $k\in K$ one has $\lexp{(g,k)}{X*Y}= \lexp{(g,h)}{X}*\lexp{(h,k)}{Y}$.

\smallskip
{\rm (b)} For a twisted diagonal subgroup $\Delta(P,\alpha,Q)$ of $G\times H$ and $(g,h)\in G\times H$, one has
$\lexp{(g,h)}{\Delta(P,\alpha,Q)}=\Delta(\lexp{g}{P}, c_g\circ\alpha\circ c_h^{-1}, \lexp{h}{Q})$.

\smallskip
{\rm (c)} For twisted diagonal subgroups $\Delta(P,\alpha,Q)$ of $G\times H$ and $\Delta(R,\beta,S)$ of $H\times K$, one has
\begin{equation*}
   \Delta(P,\alpha,Q) * \Delta(R,\beta,S) = 
   \Delta \Bigl(\alpha(Q\cap R), \alpha|_{Q\cap R}\circ \beta|_{\beta^{-1}(Q\cap R)}, \beta^{-1}(Q\cap R)\Bigr)\,.
\end{equation*}
\end{lemma}

In the following theorem, see \cite{bouc2010inducedbimodules}, we interpret any $(RG,RH)$-module $M$ also as left $R[G\times H]$-module and vice-versa via $(g,h)m= gmh^{-1}$, for $g\in G$, $h\in H$, and $m\in M$.

\begin{theorem}[Bouc]\label{thm bouc}
Let $X\le G\times H$ and $Y\le H\times K$, and let $M$ be an $RX$-module and $N$ an $RY$-module. Then, one has an isomorphism
\begin{equation*}
   \Ind_X^{G\times H}(M) \otimes_{RH} \Ind_Y^{H\times K}(N) \cong
   \bigoplus_{t\in [p_2(X)\backslash H/p_1(Y)]} \Ind_{X * \lexp{(t,1)}{Y}}^{G\times K} (M * \lexp{(t,1)}{N})
\end{equation*}
of $R[G\times K]$-modules, where $[p_2(X)\backslash H/p_1(Y)]$ denotes a set of representatives of the $(p_2(X),p_1(Y))$-double cosets in $H$.
\end{theorem}

We point out that if one replaces in the above theorem a representative $t$ by another representative of the same double coset then the resulting summands are isomorphic as $R[G\times K]$-modules.


\section{The case $G=D\rtimes E$ with $D$ cyclic and $E\le \Aut(D)$}\label{sec special G}

In this section we fix the following situation: Assume that $G=D\rtimes E$, were $D$ is a cyclic group of order $p^n$ for $n\ge 1$ and $E\le \Aut(D)$ is a $p'$-subgroup. Since $\Aut(D)$ is cyclic for odd $p$ and since $|\Aut(D)|=p^{n-1}(p-1)$ this implies that $E$ is a cyclic group of order dividing $p-1$ and that $E$ is trivial if $p=2$. We will write $D_i$ for the subgroup of $D$ of order $p^i$, for $i=0,\ldots,n$. Since $C_G(D)=D$, $1_{FG}$ is primitive in $Z(FG)$ and $FG$ is a block algebra (see Theorem~2.8 in \cite[Chapter~5]{nagao2014representations} and recall that $1_{FD}$ is primitive in $Z(FD)=FD$, since $D$ is a $p$-group). Note that the Cartan matrix of $FG$ is of the form (\ref{eqn special C}) with $m=\frac{|D|-1}{|E|}$, see \cite[Theorem~17.2]{alperin1986local}. 

We focus for the rest of this section on $T^\Delta(FG,FG)$. We write elements of $D\rtimes E$ as $x\rho$ with $x\in D$ and $\rho\in E$ with multiplication given by $(x\rho)(y\sigma)= x\rho(y)\rho\sigma$. The twisted diagonal $p$-subgroups of $G\times G$ are of the form $\Delta(D_i,\alpha,D_i)$, with $i\in\{0,\ldots,n\}$ and $\alpha\in\Aut(D_i)$. We set $\hat{E}:=\Hom(E, F^\times)$ and note that, up to isomorphism, the simple $FG$-modules are given by $F_\lambda$, $\lambda\in\hat{E}$, where $F_\lambda=F$ as $F$-vector space and $x\rho\cdot a:=\lambda(\rho)a$ for $x\rho\in G$ and $a\in F_\lambda$. By abuse of notation we also view $F_\lambda$ via restriction as a simple $FE$-module.

Let $i\in\{1,\ldots,n\}$. Since the canonical homomorphism $(\ZZ/p^n\ZZ)^\times \to (\ZZ/p^i \ZZ)^\times$ is surjective, also the restriction map $\pi_i\colon\Aut(D)\mapsto\Aut(D_i)$, $\alpha\mapsto \alpha|_{D_i}$, is surjective. Since $|\Aut(D)|=p^{n-1}(p-1)$ and $|\Aut(D_i)|=p^{i-1}(p-1)$, we have $|\ker(\pi)|=p^{n-i}$. Since $E\le \Aut(D)$ is a $p'$-subgroup, $\pi_i|_E\colon E\to\Aut(D_i)$ is injective. This implies that for any $\rho\in E$ and any $x\in D$ with $\rho(x)=x$, one has $x=1$ or $\rho=1$. In other words, $E$ acts Frobeniusly on $D$.

\smallskip
The following Lemma provides some group theoretic properties of $G\times G$ that will be needed later.

\begin{lemma}\label{lem special G}
Let $i,j\in\{0,\ldots,n\}$, $\alpha\in\Aut(D_i)$, $\beta\in\Aut(D_j)$, and $g\in G$. Further, set $k:=\min\{i,j\}$, $l:=\max\{i,j\}$. 

\smallskip
{\rm (a)}  The $p$-subgroups $\Delta(D_i,\alpha,D_i)$ and $\Delta(D_j,\beta,D_j)$ are $G\times G$-conjugate if and only if $i=j$ and there exists $\rho\in E$ with $\alpha=\beta\circ\rho|_{D_i}$.

\smallskip
{\rm (b)} Suppose that $i\ge 1$. Then one has $N_{G\times G}(\Delta(D_i,\alpha,D_i))=(D\times D)\Delta(E)$. In particular, $\Delta(D_iE,\tilde\alpha,D_iE)=\Delta(D_i,\alpha,D_i)\Delta(E)$ is a subgroup of $G\times G$, where $\tilde\alpha(x\rho):=\alpha(x)\rho$ for $x\in D_i$ and $\rho\in E$, and $\tilde\alpha$ is an automorphism of $D_iE$.

\smallskip
{\rm (c)} Suppose that $(i,j)\neq (0,0)$ and $g\notin D_lE$. Then $D_iE\cap \lexp{g}{(D_jE)}=D_k$. Moreover, $|D_iEgD_jE|=p^l |E|^2$ and the number of nontrivial double cosets in $D_iE\backslash G/D_jE$ is equal to $(p^{n-l}-1)/|E|$.
%
%
%
\end{lemma}

\begin{proof}
(a) Using the formula in Lemma~\ref{lem omnibus}(b), that $D\times D$ centralizes $\Delta(D_i,\alpha,D_i)$ and that $\Aut(D_i)$ is abelian, we immediately obtain the result.

\smallskip
(b) $D\times D$ centralizes $\Delta(D_i,\alpha,D_i)$ and since $\Aut(D_i)$ is abelian, a quick computation shows that $\Delta(E)$ normalizes $\Delta(D_i,\alpha,D_i)$. This shows one inclusion of the statement. Conversely, by Part(a) and since the restriction map $\pi_i\colon E\to\Aut(D_i)$ is injective, the $G\times G$-conjugacy class of $\Delta(D_i,\alpha,D_i)$ has at least $|E|$ elements. Since $(D\times D)\Delta(E)$ has index $|E|$ in $G\times G$, the result follows.

\smallskip
(c) We may assume without loss of generality that $i\le j$ by noting that $D_iE \cap \lexp{g}{(D_jE)} = D_k$ if and only if $\lexp{g^{-1}}{(D_iE)}\cap D_jE= D_k$ and by noting that taking inverses induces a bijection between $D_iEgD_jE$ and $D_jEg^{-1}D_iE$ and between $D_iE\backslash G/ D_jE$ and $D_jE\backslash G/D_iE$ fixing the trivial double coset $D_iED_jE=D_iD_jE=D_lE$. 

To prove the first statement, write $g=z \rho$ with $z\in D$ and $\rho\in E$. Since $\lexp{g}(D_jE)= D_j\cdot\lexp{g}{E}=D_j\lexp{z}{E}$, we may also assume that $g=z\in D$, but $z\notin D_j$. 
Then, clearly $D_i\le D_iE\cap D_j \lexp{z}{E}$. 
Conversely consider an element $x\rho=y\lexp{z}{\sigma}\in D_iE\cap D_j\lexp{z}{E}$ with $x\in D_i$, $y\in D_j$ and $\rho,\sigma\in E$. Since $y\lexp{z}{\sigma} = yz\sigma(z^{-1})\sigma$ we obtain $x=yz\sigma(z^{-1})$ and $\rho=\sigma$. 
Since $\sigma|_{D_1}\in\Aut(D_1)$, $\sigma(z)=z^r$ for some $r\in\ZZ$ with $p\nmid r$. 
Thus, we obtain $x=yz^{1-r}$. Since $xy^{-1}\in D_j$ but $z\notin D_j$, we obtain $r\equiv 1\mod p$ and therefore $\sigma|_{D_1}=\id_{D_1}$. 
The injectivity of the restriction map $\pi_1\colon E\to\Aut(D_1)$ now implies $\rho=\sigma=1$ and $x\rho=x\in D_i$ as claimed.

For the second statement let $g\in G$ with $g\notin D_iED_jE=D_iD_jE=D_jE$. Then, using the first statement, we obtain 
\begin{equation*}
   |D_iEgD_jE| = |D_iE\lexp{g}{(D_jE)}|=\frac{|D_iE|\cdot|D_jE|}{|D_iE\cap \lexp{g}{(D_jE)}|} = \frac{|D_iE|\cdot|D_jE|}{|D_i|} =
   p^{j}|E|^2\,.  
\end{equation*}
Therefore the $p^n|E|-p^j|E|$ elements in $G\smallsetminus D_jE$ partition into the non-trivial $(D_iE,D_jE)$-double cosets of constant size $p^j|E|^2$ and the second statement follows, since $i=k$ and $j=l$.
%
%
%
\end{proof}

The following two propositions provide a $\ZZ$-basis of $T^\Delta(FG,FG)$ and a formula for the product of any two basis elements. By Theorem~\ref{thm normal Sylow}, for $\lambda\in\hat E$, the $FG$-module $P_\lambda:=\Ind_E^G(F_\lambda)$ is the projective cover of the simple $FG$-module $F_\lambda$. Thus, setting $P_{\lambda,\mu} := P_\lambda\otimes_F P_{\mu^{-1}} \cong \Ind_{E\times E}^{G\times G}(F_\lambda\otimes F_{\mu^{-1}})$, for $\lambda,\mu\in\hat{E}$, is consistent with the notation introduced in \ref{noth Pij}.

\begin{proposition}\label{prop basis}
{\rm (a)} The projective indecomposable $F[G\times G]$-modules are of the form $P_{\lambda,\mu}$ for $\lambda,\mu\in\hat{E}$. Moreover, $P_{\lambda,\mu}\cong P_{\lambda',\mu'}$ if and only if $\lambda=\lambda'$ and $\mu=\mu'$.

\smallskip
{\rm (b)} Let $i\in\{1,\ldots,n\}$, $\alpha\in\Aut(D_i)$ and $\lambda\in \hat{E}$. Then $M_{i,\alpha,\lambda}:=\Ind_{\Delta(D_iE,\tilde\alpha, D_iE)}^{G\times G}(F_\lambda)$ is an indecomposable $F[G\times G]$-module with vertex $\Delta(D_i,\alpha,D_i)$, where by abuse of notation we view $F_\lambda$ as $F\Delta(D_iE,\tilde\alpha,D_iE)$-module via the isomorphism $E\cong \Delta(D_iE,\tilde\alpha,D_iE)/\Delta(D_i,\alpha,D_i)$, $\rho\mapsto (\rho,\rho)\Delta(D_i,\alpha,D_i)$. Moreover, every indecomposable non-projective $F[G\times G]$-module with twisted diagonal vertex is of this form, and $M_{i,\alpha,\mu}\cong M_{j,\beta,\lambda}$ if and only if $i=j$, $\alpha=\beta\circ \rho|_{D_i}$ for some $\rho\in E$, and $\lambda=\mu$.
\end{proposition}

\begin{proof}
Part~(a) is immediate and
Part~(b) follows from Theorem~\ref{thm normal Sylow} together with Lemma~\ref{lem special G}(b) and noting that $\Delta(E)$ is a complement of $D\times D=N_{D\times D}(\Delta(D_i,\alpha,D_i))$ in $N_{G\times G}(\Delta(D_i,\alpha,D_i))=(D\times D)\Delta(E)$.
\end{proof}

\begin{proposition}\label{prop mult formulas}
Let $\lambda,\mu,\lambda',\mu'\in\hat E$, $i,j\in\{1,\ldots,n\}$, $\alpha\in\Aut(D_i)$, and $\beta\in\Aut(D_j)$. Then the standard basis elements of $T^\Delta(FG,FG)$ satisfy the following multiplication rules:
\begin{align*}
   [P_{\lambda,\mu}]\cdot_G [P_{\lambda',\mu'}] = & \,
            \begin{cases} (m+1) [P_{\lambda,\mu'}] & \text{if $\mu=\lambda'$,} \\ m [P_{\lambda,\mu'}] & \text{if $\mu\not=\lambda'$,}\end{cases} \\
   [P_{\lambda,\mu}]\cdot_G [M_{j,\beta,\mu'}] = & \, 
            [P_{\lambda,\mu {\mu'}^{-1}}]\, + \, \sum_{\nu\in\hat E} \frac{p^{n-j}-1}{|E|} [P_{\lambda,\nu}]\,, \\
   [M_{i,\alpha,\lambda}] \cdot_G [P_{\lambda',\mu'}] = & \, 
            [P_{\lambda\lambda',\mu'}]\, + \, \sum_{\nu\in\hat E} \frac{p^{n-i}-1}{|E|} [P_{\nu,\mu'}]\,,  \\
   [M_{i,\alpha,\lambda}] \cdot_G [M_{j,\beta,\mu}] = & \, 
           [M_{k,\alpha|_{D_k}\circ\beta|_{D_k},\lambda\mu}]\, + \,
      \sum_{\nu\in\hat E} \frac{p^{n-l}-1}{|E|} [M_{k, \alpha|_{D_k}\circ\beta|_{D_k}, \nu}] \,,  \\
\end{align*}
where $m=(|D|-1)/|E|$, $k=\min\{i,j\}$, and $l=\max\{i,j\}$. Moreover, $[M_{n,\id,1}]$ is the identity element of $T^\Delta(FG,FG)$.
\end{proposition}

\begin{proof}
The first equation follows from Lemma~\ref{lem Cartan}(a).

To prove the second equation we use Theorem~\ref{thm bouc} and obtain
\begin{align*}
   \Ind_{E\times E}^{G\times G} & 
   (F_\lambda\otimes_F F_{\mu^{-1}}) \otimes_{FG} \Ind_{\Delta(D_jE, \tilde\beta, D_jE)}^{G\times G} (F_{\mu'}) \cong \\
   &   \bigoplus_{t \in [E\backslash G/D_jE]} \Ind_{(E\times E)\, * \lexp{(t,1)}{\Delta(D_jE,\betatilde,D_jE)}}^{G\times G} 
   \bigl((F_\lambda\otimes_F F_{\mu^{-1}}) * \lexp{(t,1)}{F_{\mu'}}\bigr)
\end{align*}
For the trivial double coset $ED_jE=D_jE$, we choose the representative $t=1$. A quick computation shows that $(E\times E) * \Delta(D_jE,\betatilde,D_jE) = E\times E$. Moreover, $(\rho,\sigma)\in E\times E$, with connecting element $\sigma$, acts on $ (F_\lambda\otimes_F F_{\mu^{-1}}) \otimes_F F_{\mu'}$ by multiplication with $\lambda(\rho)\mu^{-1}(\sigma) \mu'(\sigma)$. Thus, $(F_\lambda\otimes_F F_{\mu^{-1}}) * F_{\mu'}\cong F_{\lambda}\otimes_F F_{\mu^{-1}\mu'}$ as $F[E\times E]$-modules. Therefore, the trivial double coset contributes the summand $P_{\lambda,\mu {\mu'}^{-1}}$. We may choose the nontrival double coset representatives as elements $t\in D$, with $t\notin D_i$. Lemma~\ref{lem special G}(c) implies that $(E\times E) * \lexp{(t,1)}{\Delta(D_jE,\betatilde,D_jE)} = E\times \{1\}$, since $E\cap\lexp{t}{D_jE}=\{1\}$. Moreover, $(\rho,1)\in E\times \{1\}$, with connecting element $1$, acts on $(F_\lambda\otimes_F F_{\mu^{-1}}) \otimes_F \lexp{(t,1)}{F_{\mu'}}$ via multiplication with $\lambda(\rho)$. Thus, $(F_\lambda\otimes_F F_{\mu^{-1}}) * \lexp{(t,1)}{F_{\mu'}}\cong F_\lambda\otimes_F F$ as $F[E\times \{1\}]$-modules. Induction to $E\times E$ results in the $F[E\times E]$-module $\bigoplus_{\nu\in\hat E} F_\lambda\otimes_F F_\nu$, and further induction to $G\times G$ yields $\bigoplus_{\nu\in\hat E} P_{\lambda,\nu}$, which is independent of $t$. Since there are precisely $(p^{n-j}-1)/|E|$ nontrivial double cosets, see Lemma~\ref{lem special G}(c), the result follows.

The third equation is proved in a similar way as the second.

For the fourth equation, Theorem~\ref{thm bouc} yields
\begin{align*}
    M_{i,\alpha,\lambda} \otimes_{FG} M_{j,\beta,\mu} \cong & \
       \Ind_{\Delta(D_iE,\tilde\alpha,D_iE)}^{G\times G} (F_\lambda) \otimes_{FG} 
       \Ind_{\Delta(D_jE,\tilde\beta,D_jE)}^{G\times G}(F_\mu) \\
    \cong & \bigoplus_{t\in [D_iE\backslash G/D_jE]}
       \Ind_{\Delta(D_iE,\tilde\alpha,D_iE)\, * \lexp{(t,1)}{\Delta(D_jE,\tilde\beta,D_jE)}}^{G\times G} (F_\lambda * \lexp{(t,1)}{F_\mu})\,.
\end{align*}
For the trivial double coset in $D_iE\backslash G/D_jE$ we choose the representative $t=1$. 
By Lemma~\ref{lem omnibus}(c), we have $\Delta(D_iE,\tilde\alpha,D_iE) * \Delta(D_jE,\tilde\beta,D_jE) = \Delta(D_kE,\tilde\alpha|_{D_kE}\circ \tilde\beta|_{D_kE}, D_kE)$. 
Moreover, the element $(\alpha(\beta(x))\rho,x\rho)\in \Delta(D_kE,\tilde\alpha|_{D_kE}\circ \tilde\beta|_{D_kE}, D_kE)$, with connecting element $\beta(x)\rho$, acts on $F_\lambda * F_\mu =F_\lambda\otimes_F F_\mu$ by multiplication with  $\lambda(\rho)\mu(\rho)$. 
Thus, $F_\lambda * F_\mu \cong F_{\lambda\mu}$ as $F\Delta(D_kE,,\tilde\alpha|_{D_kE}\circ \tilde\beta|_{D_kE}, D_kE)$-modules and the trivial double coset contributes the direct summand $M_{k,\alpha|_{D_k}\circ\beta|_{D_k}, \lambda\mu}$. 
We may choose the nontrivial double coset representatives as elements $t\in D$ with $t\notin D_l$. 
Lemma~\ref{lem special G}(f) implies $\Delta(D_iE,\tilde\alpha,D_iE) * \lexp{(t,1)}{\Delta(D_jE,\tilde\beta,D_jE)} = \Delta(D_k,\alpha|_{D_k}\circ\beta_{D_k}, D_k)$. Moreover, the element $(\alpha(\beta(x)), x)\in \Delta(D_k,\alpha|_{D_k}\circ\beta_{D_k}, D_k)$, with connecting element $\beta(x)$, acts on $F_\lambda * \lexp{(t,1)}{F_\mu} = F_\lambda\otimes_F \lexp{(t,1)}{F_\mu}  $ by multiplication with $1$, i.e., as the identity. 
Thus, $F_\lambda * \lexp{(t,1)}{F_\mu}$ is the trivial $F\Delta(D_k,\alpha|_{D_k}\circ\beta|_{D_k}, D_k)$-module and induction to $\Delta(D_kE, \tilde\alpha|_{D_kE}\circ\tilde\beta|_{D_kE}, D_kE)$ yields $\bigoplus_{\nu\in \hat E} F_\nu$. Further induction to $G\times G$ yields $\bigoplus_{\nu\in\hat E} M_{k,\alpha|_{D_k}\circ\beta|_{D_k}, \nu}$, independent of $t$. Since there are precisely $(p^{n-l}-1)/|E|$ nontrivial double cosets in $D_iE\backslash G/D_jE$, see Lemma~\ref{lem special G}(c), the result follows.

\smallskip
For the proof of the  last statement, note that $M_{n,\id,1}=\Ind_{\Delta(G)}^{G\times G}(F)\cong FG$ as $(FG,FG)$-bimodules, or alternatively use the multiplication formulas.
\end{proof}

Recall the definition of $T_i^\Delta(FG,FG)$ from Section~\ref{sec T(FG)}. The following corollary follows immediately from Proposition~\ref{prop mult formulas}.

\begin{corollary}\label{cor mult formulas}
For $i\in\{0,\ldots,n\}$, the subgroup $T^\Delta_i(FH,FH)$ of $T^\Delta(FH,FH)$ is multiplicatively closed and the ring $T^\Delta(B,B)/Pr(B,B)$ is commutative. More precisely, the following hold:

\smallskip
{\rm (a)} $[P_{\lambda,\mu}]\mapsto E_{\lambda,\mu}$ defines a ring isomorphism $Pr(FH,FH)\cong\Mat_{|E|}(\ZZ)_C$, where $C$ is the matrix in (\ref{eqn special C}) with $m=(|D|-1)/|E|$ and $l=|E|$;

\smallskip
{\rm (b)} The map $[M_{n,\alpha,\lambda}]\mapsto (\alpha E, \lambda)$ induces a ring isomorphism $T^\Delta_n(FG,FG)\cong \ZZ[\Aut(D)/E \times \hat E]$;

\smallskip
{\rm (c)} For $i\in \{1,\ldots,n-1\}$, the map $[M_{i,\alpha,\lambda}]\mapsto (\alpha\pi_i(E),\lambda)$ induces a ring isomorphism $T^\Delta_i(FG,FG)\cong \ZZ[\Aut(D_i)/\pi_i(E)\times \hat E]_{c_i}$, where $c_i=(1,1) + m_i \sum_{\lambda\in\hat E} (1,\lambda)$ with $m_i=(p^{n-i}-1)/|E|$;
\end{corollary}

Before we take the general study of idempotents in $T^\Delta(B,B)$ further, we will establish two general ring-theoretic results.


\section{Two ring theoretic lemmas}\label{sec gen idem}

In this section, we prove two ring-theoretic lemmas that will be used in the proofs of Theorem~C and D in Section~\ref{sec proof of Thm C}.

\smallskip
Note that if $R$ is a unitary Noetherian ring, then every idempotent $e$ of $R$ has a primitive decomposition.

\begin{lemma}\label{lem chain of ideals}
Let $R$ be a unital Noetherian ring, let $n\ge 0$, and suppose that $I_0\subseteq I_1\subseteq\ldots\subseteq I_n\subset R$ is a chain of ideals of $R$ such that $\idem(R/I_n)=\{0,1\}$ and $\idem(I_k/I_{k-1})=\{0\}$ for all $k=1,\ldots,n$.

\smallskip
{\rm (a)} For any two nonzero orthogonal idempotents $e,e'$ of $R$ one has $e\in I_0$ or $e'\in I_0$. In particular, $\idem(R/I_0)=\{0,1\}$.

\smallskip
{\rm (b)} If $\{e_1,\ldots, e_r\}\subseteq \idem(I_0)$ is a maximal set of pairwise orthogonal primitive idempotents of $I_0$ then $\{e_1,\ldots,e_r, 1-(e_1+\cdots+e_r)\}$ is a primitive decomposition of $1_R$ in $R$.
\end{lemma}

\begin{proof}
Induction on $n$. Assume first that $n=0$ and let $e,e'\in \idem(R)$ be non-zero and orthogonal. Then $e+I_0$ and $e'+I_0$ are orthogonal in $\idem(R/I_0)=\{0,1\}$.  Therefore, $e\in I_0$ or $e'\in I_0$, establishing (a). To see (b), it suffices to show that the idempotent $1-e$ is primitive in $R$, where $e:=e_1+\cdots +e_r$. So suppose that $1-e=f_1+f_2$ with orthogonal non-zero idempotents $f_1,f_2$ of $R$. By (a), $f_1$ or $f_2$ is in $I_0$. By symmetry we may assume that $f_1\in I_0$. Let $f$ be an element of a primitive decomposition of $f_1$ in $R$. With $f_1$ also $f$ is in $I_0$ and orthogonal to $e$. Thus, $e_1,\ldots,e_r,f$ are primitive pairwise orthogonal idempotents of $I_0$, contradicting the maximality of $\{e_1,\ldots,e_r\}$. 

Now suppose that $n>0$ and consider the chain of ideals $I_n/I_{n-1}\subset R/I_{n-1}$ of $R/I_{n-1}$ which satisfies the hypothesis of the lemma for the base case. Since the empty set is a maximal set of pairwise orthogonal primitive idempotents of $I_n/I_{n-1}$, we obtain that $1_R$ is primitive in $R$, i.e., $\idem(R/I_{n-1})=\{0,1\}$. Therefore, also the shorter chain $I_0\subseteq I_1\subseteq\cdots\subseteq I_{n-1}\subset R$ satisfies the hypothesis of the lemma and the result follows by induction. 
\end{proof}

In the proof of the next lemma we will need the following result. Recall the definition of $R_c$ in \ref{not Rc}.

\begin{theorem}(Coleman, \cite{coleman1966idempotents})\label{thm coleman}
Let $G$ be a finite group and $R$ an integral domain in which no prime divisor of $|G|$ is invertible. Then $\idem(R[G])=\{0,1\}$.
\end{theorem}

\begin{lemma}\label{lem only 0 idempotent}
Let $\Gamma$ and $X$ be finite groups, $R$ a commutative ring, $m\in R$, and $c:=(1,1)+m\sum_{x\in X} (1,x) \in R[\Gamma\times X]$.

\smallskip
{\rm (a)} Suppose that $R$ is an integral domain such that no prime divisor of $|\Gamma\times X|$ is invertible in $R$ and $1+m|X|\notin R^\times\cup\{0\}$. Then $\idem(R[\Gamma\times X]_c)=\{0\}$.

\smallskip
{\rm (b)} Suppose that $1+m|X|$ is invertible in $R$. Then $d:=(1,1) - m(1+m|X|)^{-1}\sum_{x\in X}(1,x)\in R[\Gamma\times X]$ is an inverse of $c$ in $R[\Gamma\times X]$. In particular, $d$ is an identity element of $R[\Gamma\times X]_c$.
\end{lemma}

\begin{proof}
(a) It is straightforward to check that the $R$-linear extensions of the maps $\Gamma\times X\to R[\Gamma\times \{1\}]$, $(\gamma,x)\mapsto (1+|X|m)(\gamma,1)$ and $(\gamma,x)\mapsto (\gamma,1)$ define $R$-linear ring homomorphism $\epsilon'\colon R[\Gamma\times X]_c\to R[\Gamma\times \{1\}]$ and $\epsilon\colon R[\Gamma\times X]\to R[\Gamma\times \{1\}]$, respectively. It is also straightforward to check that 
\begin{equation}\label{eqn * mult}
   a\cdot_cb= a\cdot b+m\epsilon(a)\cdot\epsilon(b)\cdot s\quad\text{for all $a,b\in R[\Gamma\times X]$, where $s:=\sum_{x\in X}(1,x)$.}
\end{equation}
In fact, since both sides are $R$-bilinear in $a$ and $b$, it suffices to check this for the standard $R$-basis elements.
Let $e$ be an idempotent of $R[\Gamma\times X]_c$. Then $\epsilon'(e)$ is an idempotent in $\epsilon'(R[\Gamma\times X]_c)= (1+m|X|)R[\Gamma\times \{1\}]$. By Theorem~\ref{thm coleman}, $\idem(R[\Gamma])=\{0,1\}$. Since $(1+m|X|)$ is not a unit in $R$, this implies that $\epsilon'(e)=0$. Since $R[\Gamma\times\{1\}]$ is $R$-torsionfree and $1+m|X]\neq 0$ we also obtain $\epsilon(e)=0$.
Equation~(\ref{eqn * mult}) now implies that $e = e \cdot_c e = e\cdot e +m\epsilon(e)\cdot\epsilon(e)\cdot s = e \cdot e$ is an idempotent also in the ring $R[\Gamma\times X]$.
Again by Theorem~\ref{thm coleman}, this implies that $e\in \{0,1\}$. Since $\epsilon(e)=0$, this implies that $e\neq 1$ and the proof is complete.

\smallskip
(b) This is a straightforward computation.
\end{proof}


\section{Proofs of Theorem~C and D}\label{sec proof of Thm C}

Throughout this section we assume that $B$ is a block algebra of $FG$ with cyclic defect group $D$ of order $p^n>1$. If $D$ is trivial then $T^\Delta(B,B)=Pr(B,B)\cong \ZZ$ and Theorems~C and D hold trivially. We choose a maximal $B$-Brauer pair $(D,e)$ and set $E:=N_G(D,e)/C_G(D)$. Then, via the conjugation action of $N_G(D,e)$ on $D$, the group $E$ can be considered as a subgroup of $\Aut(D)$. It is known that $E$ is a $p'$-group, see \cite[Theorem~6.7.6(v)]{Li18b}. We adopt all the notation introduced at the beginning of Section~\ref{sec special G}, except, that in this section we set $H:=D\rtimes E$.

\smallskip
Many questions concerning blocks with cyclic defect groups reduce to the group $H$. This also holds for the ring structure of $T^\Delta(B,B)$.

\begin{theorem}\label{thm reduction}
There exists a ring isomorphism $T^\Delta(B,B)\cong T^\Delta(FH,FH)$ that restricts to an isomorphism $T^\Delta_{\le i}(B,B)\cong T^\Delta_{\le i}(FH,FH)$ for all $i\in \{0,\ldots,n\}$.
\end{theorem}

\begin{proof}
By \cite{rouquier1998cyclic}, there exists a splendid Rickard equivalence between the block $B$ and its Brauer correspondent block $B_1\in\Bl(FN_G(D))$. By \cite[Theorem~1.5]{boltjexu2008p-permutation} this implies the existence of a $p$-permutation equivalence $\gamma\in T^\Delta(B,B_1)$. It is straightforward to verify that the map $\gamma^\circ\cdot_G - \cdot_G \gamma\colon T^\Delta(B,B)\to T^\Delta(B_1,B_1)$ is a ring isomorphism. 

Furthermore, if $i\in (F[N_G(D)])^D$ is a source idempotent of $B_1$, then $iB_1i$ is isomorphic to $FH$ as interior $D$-algebra, see \cite[Theorem~6.14.1]{Li18a}. Note that the $2$-cocycle appearing there is trivial, since $E$ is cyclic. One can consider $iB_1$ as $(FH,B_1)$-bimodule and it induces a Morita equivalence between $B_1$ and $FH$. Moreover, viewed after restriction as $F[D\times N_G(D)]$-module, $iB_1$ is a direct summand of the permutation module $FN_G(D)$. This implies that the $(F[H\times N_G(D)])$-module $iB_1$ is a $p$-permutation module. It follows that one obtains a ring isomorphism $T^\Delta(B_1,B_1)\cong T^\Delta(FH,FH)$ as in the previous paragraph.

Theorem~\ref{thm bouc} implies that the above bijections map $T^\Delta_{\le i}(B,B)$ to $T^\Delta_{\le i}(FH,FH)$ and vice-versa.
\end{proof}

\begin{proof} {\em of Theorem~C}.\quad
By Theorem~\ref{thm reduction} it suffices to prove the statement in the case that $G=H=D\rtimes E$ as introduced at the beginning of the section. In this case we consider the chain of ideals
\begin{equation*}
  \{0\}\subset T_0^\Delta(FG,FG) \subset T_{\le 1}^\Delta(FG,FG) \subset \cdots \subset T_{\le n-1}^\Delta(FG,FG) \subset T_{\le n}^\Delta(FG,FG)=T^\Delta(FG,FG)\,.
\end{equation*}
By Corollary~\ref{cor mult formulas}, the ring $T^\Delta(FG,FG)/Pr(FG,FG)$ is abelian. By Example~\ref{ex special C}, the multiplicity of $1$ as elementary divisor of the Cartan matrix of $FG$ is equal to $l-1=|E|-1$. 
Therefore, using Theorem~\ref{thm Pr(B,B)}(b) and Lemma~\ref{lem chain of ideals}, for the remaining statements of Theorem~C, it suffices to show that

\smallskip
(i) $\idem\bigl(T^\Delta(FG,FG)/T^\Delta(FG,FG)_{\le n-1}\bigr) = \{0,1\}$ and

\smallskip
(ii) $\idem\bigl(T^\Delta(FG,FG)_{\le i}/T^\Delta(FG,FG)_{\le i-1} \bigr) = \{0\}$ for all $i\in \{1,\ldots,n-1\}$.

\smallskip
But (i) follows immediately from Corollary~\ref{cor mult formulas}(b) and Theorem~\ref{thm coleman}, and (ii) follows immediately from Corollary~\ref{cor mult formulas}(c) and Lemma~\ref{lem only 0 idempotent}(a).
\end{proof}

\begin{proof} {\em of Theorem~D}.\quad 
Again, by Theorem~\ref{thm reduction} it suffices to prove Theorem~D in the case that $G=H=D\rtimes E$ and from now on we assume this. Furthermore we assume that the characteristic of $k$ is not $p$.

\smallskip
{\em Claim 1: For each $i\in\{0,\ldots,n\}$, $kT^\Delta_i(FG,FG)$ has an identity element $e_i$.}\quad
We define
\begin{equation*}
   e_0:=\sum_{\lambda\in \hat E} [P_{\lambda,\lambda}] - \frac{m}{|D|} \sum_{\lambda,\mu\in \hat E} [P_{\lambda,\mu}]\,.
\end{equation*}
Then, by Example~\ref{ex special C}, $e_0$ is an identity element of $kT_0^\Delta(FG,FG)=k Pr(FG,FG)$. Further, for $i\in \{1,\ldots,n\}$, we set 
\begin{gather*}
   m_i:=(p^{n-i}-1)/|E|\,,\ m'_i:= - m_i(1+m_i|E|)^{-1} = -m_i/ p^{n-i}\in k \quad\text{and}  \\
   c_i:=[M_{i,\id,1}] + m_i\sum_{\lambda\in \hat E} [M_{i,\id,\lambda}]\,,\
   e_i:=[M_{i,\id,1}] + m'_i\sum_{\lambda\in \hat E} [M_{i,\id,\lambda}]\in kT_i^\Delta(FG,FG)\,.
\end{gather*}
Note that $m_n=m'_n=0$. Using the isomorphisms in Corollary~\ref{cor mult formulas}(b) and (c) together with Lemma~\ref{lem only 0 idempotent}(b), we see that $e_i$ is an identity element of $kT_i^\Delta(FG,FG)$. Claim~1 is now established.

\smallskip
{\em Claim 2: For any $i\in\{0,\ldots,n\}$, the idempotent $e_i$ is central in $T^\Delta(FG,FG)$.}\quad To see that $e_0$ is central, it suffices to show $[M_{i,\alpha,\lambda}]\cdot e_0= e_0 \cdot [M_{i,\alpha,\lambda}]$ for any $i\in \{1,\ldots,n\}$, $\alpha\in\Aut(D_i)$, and $\lambda\in \hat E$. We compute both sides, using Proposition~\ref{prop mult formulas}, and observe that they coincide:
\begin{align*}
   [M_{i,\alpha,\lambda}]\cdot e_0 & = \sum_{\mu\in\hat E} [M_{i,\alpha,\lambda}] \cdot [P_{\mu,\mu}] - 
                \frac{m}{|D|} \sum_{\mu,\nu\in\hat E} [M_{i,\alpha,\lambda}] \cdot [P_{\mu,\nu}] \\
   & = \sum_{\mu\in \hat E} \Bigl([P_{\lambda\mu,\mu}] + m_i\sum_{\nu\in\hat E} [P_{\nu,\mu}]\Bigr) -
            \frac{m}{|D|} \sum_{\mu,\nu\in\hat E} \Bigl( [P_{\lambda\mu, \nu}] + m_i\sum_{\nu'\in \hat E} [P_{\nu',\nu}] \Bigr) \\
   & = \sum_{\mu\in\hat E} [P_{\lambda\mu,\mu}] + \Bigl(m_i-\frac{m}{|D|} -\frac{mm_i|E|}{|D|}\Bigr) \sum_{\mu,\nu\in\hat E} [P_{\mu,\nu}]
\end{align*}
and
\begin{align*}
   e_0\cdot [M_{i,\alpha,\lambda}] & = \sum_{\mu\in\hat E} [P_{\mu,\mu}]\cdot[M_{i,\alpha,\lambda}] - 
         \frac{m}{|D|} \sum_{\mu,\nu} [P_{\mu,\nu}] \cdot [M_{i,\alpha,\lambda}] \\
   & = \sum_{\mu\in\hat E} \Bigl( [P_{\mu,\mu\lambda^{-1}}] + m_i \sum_{\nu\in\hat E} [P_{\mu,\nu}] \Bigr) -
         \frac{m}{|D|} \sum_{\mu,\nu\in\hat E} \Bigl( [P_{\mu,\nu\lambda^{-1}}] + m_i \sum_{\nu'\in\hat E} [P_{\mu,\nu'}] \Bigr) \\
   & = \sum_{\mu\in\hat E} [P_{\mu,\mu\lambda^{-1}}] + \Bigl( m_i - \frac{m}{|D|} - \frac{mm_i|E|}{|D|} \Bigr) 
            \sum_{\mu,\nu\in \hat E} [P_{\mu,\nu}]\,.
\end{align*}
Now let $i\in\{1,\ldots,n\}$. Since any two elements in $\bigoplus_{i=1}^n T_i^\Delta(FG,FG)$ commute by the fourth equation in Proposition~\ref{prop mult formulas}, it suffices to show that $[P_{\lambda,\mu}]\cdot e_i = e_i \cdot [P_{\lambda,\mu}]$ for all $\lambda,\mu\in\hat E$. Again, a straightforward computation using Proposition~\ref{prop mult formulas} yields
\begin{equation*}
   [P_{\lambda,\mu}]\cdot e_i = [P_{\lambda,\mu}] + 
   \Bigl(m_i - \frac{m_i}{p^{n-i}} - \frac{m_i^2|E|}{p^{n-i}} \Bigr) \sum_{\nu\in \hat E} [P_{\lambda,\nu}]
\end{equation*}
and
\begin{equation*}
   e_i \cdot [P_{\lambda\mu}] = [P_{\lambda,\mu}] + 
   \Bigl(m_i - \frac{m_i}{p^{n-i}} - \frac{m_i^2|E|}{p^{n-i}} \Bigr) \sum_{\nu\in\hat E} [P_{\nu,\mu}]\,.
\end{equation*}
But a quick computation shows that 
\begin{equation}\label{eqn m_i relation}
   m_i - \frac{m_i}{p^{n-i}} - \frac{m_i^2|E|}{p^{n-i}} = 0
\end{equation}
so that we obtain
\begin{equation}\label{eqn e_i central}
   [P_{\lambda,\mu}]\cdot e_i = [P_{\lambda,\mu}] = e_i \cdot [P_{\lambda\mu}] \,.
\end{equation}
   
\smallskip
{\em Claim 3: For all $i\le j$ in $\{0,\ldots,n\}$ one has $e_ie_j=e_i$.}\quad
If $i=0$ this follows from $e_0^2=e_0$ if $j=0$ and from Equation~(\ref{eqn e_i central}). If $1\le i\le j$, another straightforward computation using Proposition~\ref{prop mult formulas} yields
\begin{align*}
   e_i \cdot e_j & = \Bigl( [M_{i,\id,1}] + m'_i \sum_{\lambda\in\hat E} [M_{i,\id,\lambda}] \Bigr) \cdot
               \Bigl( [M_{j,\id,1}] + m'_i \sum_{\mu\in\hat E} [M_{j,\id,\mu}] \Bigr) \\
   & = [M_{i,\id,1}] + \Bigl( m_j + m'_i + m'_im_j|E| + m'_j + m'_jm_j|E| + m'_i m'_j|E| + m'_i m'_j m_j|E|^2 \Bigr) \sum_{\nu\in\hat E} [M_{i,\id,\nu}]\,.
\end{align*}
Using Equation~(\ref{eqn m_i relation}) for $j$ instead of $i$, we obtain $m_j + m'_j + m_jm'_j|E|=0$, which implies that the expression inside the parentheses is equal to $m'_i$. Thus, $e_i\cdot e_j = e_i$ as claimed.

\smallskip
Claims 1, 2, and 3 imply now that $e_0,e_1-e_0,\ldots, e_n-e_{n-1}$ are pairwise orthogonal central idempotents of $kT^\Delta(FG,FG)$ whose sum is equal to $1_{kT^\Delta(FG,FG)}=e_n$.

\smallskip
We already know that $kT^\Delta(FG,FG)e_0 = kPr(FG,FG) \cong \Mat_{|E|}(k)$ from Theorem~\ref{thm QPr(B,B)}. Our next goal is 

\smallskip
{\em Claim 4: For $i\in\{1,\ldots,n\}$, the $k$-algebra $kT^\Delta(FG,FG) (e_i-e_{i-1})$ is isomorphic to $k[\Aut(D_i)/\pi_i(E) \times \hat E]$.}\quad It suffices to show that the $k$-linear and multiplicative map 
\begin{equation}\label{eqn key map}
kT_i^\Delta(FG,FG)\to kT^\Delta(FG,FG) (e_i-e_{i-1})\,, \quad a\mapsto a (e_i-e_{i-1}) \,,
\end{equation}
is an isomorphism, since $kT_i^\Delta(FG,FG)$ is isomorphic to $k[\Aut(D_i)/\pi_i(E) \times \hat E]$ by Corollary~\ref{cor mult formulas}(b),(c) and Lemma~\ref{lem only 0 idempotent}(b), using that $(1+m_i|E|)=p^{n-i}$ is invertible in $k$. 
Using Proposition~\ref{prop mult formulas} and that $e_i$ is an identity element of $kT^\Delta_i(FG,FG)$, we obtain, for any $\alpha\in\Aut(D_i)$ and any $\lambda\in\hat E$, that $[M_{i,\alpha,\lambda}]\cdot(e_i-e_{i-1})  \in [M_{i,\alpha,\lambda}] + kT^\Delta_{i-1}(FG,FG)$. 
Therefore the map in (\ref{eqn key map}) followed up by the projection onto $kT^\Delta_i(FG,FG)$ is the identity map, so that the map in (\ref{eqn key map}) is injective. But then we have 
\begin{align*}
   \dim_k T^\Delta(FG,FG) & = \sum_{i=0}^n \dim_k T_i^\Delta(FG,FG) \\
   & \le \dim_k kT^\Delta(FG,FG)e_0 + \sum_{i=1}^n \dim_k kT^\Delta(FG,FG) (e_i-e_{i-1}) = \dim_k kT^\Delta(FG,FG),
\end{align*}
since $e_0, e_1-e_0,\ldots, e_n-e_{n-1}$ are pairwise orthogonal idempotents whose sum is $1$. Therefore equality must hold for each of the 
summands above. This implies that the map in (\ref{eqn key map}) is also surjective and Claim~4 is proved.

\smallskip
Finally, we prove the last statement in Theorem~D for an arbitrary field $k$. By (\ref{eqn decomposition}) and Maschke's theorem, the statement holds in the case that $k$ has characteristic different from $p$. So assume that $k$ has characteristic $p$. Since $kT^\Delta(FG,FG)/kT_{\le n-1}^\Delta(FG,FG)$ is isomorphic to $k[\Aut(D)/E\times \hat E]$ by Corollary~\ref{cor mult formulas}(b), and since $p$ divides $|\Aut(D)/E|$, the algebra $kT^\Delta(FG,FG)$ has a factor algebra that is not semisimple. Therefore, $kT^\Delta(FG,FG)$ is not semisimple and the proof of Theorem D is complete.
\end{proof}


\bibliographystyle{plain}
\bibliography{ref}

\begin{thebibliography}{10}

\bibitem{alperin1986local}
J.~L. Alperin.
\newblock {\em Local representation theory}, volume~11 of {\em Cambridge
  Studies in Advanced Mathematics}.
\newblock Cambridge University Press, Cambridge, 1986.
\newblock Modular representations as an introduction to the local
  representation theory of finite groups.

\bibitem{boltje2020ppermutation}
Robert Boltje and Philipp Perepelitsky.
\newblock $p$-permutation equivalences between blocks of group algebras, 2020.
  arXiv 2007.09253.

\bibitem{boltjexu2008p-permutation}
Robert Boltje and Bangteng Xu.
\newblock On {$p$}-permutation equivalences: between {R}ickard equivalences and
  isotypies.
\newblock {\em Trans. Amer. Math. Soc.}, 360(10):5067--5087, 2008.

\bibitem{bouc2010inducedbimodules}
Serge Bouc.
\newblock Bisets as categories and tensor product of induced bimodules.
\newblock {\em Appl. Categ. Structures}, 18(5):517--521, 2010.

\bibitem{boucyilmaz2022diagonal}
Serge Bouc and Deniz Y\i~lmaz.
\newblock Diagonal {$p$}-permutation functors, semisimplicity, and functorial
  equivalence of blocks.
\newblock {\em Adv. Math.}, 411:Paper No. 108799, 54, 2022.

\bibitem{broue1990isometries}
Michel Brou\'{e}.
\newblock Isom\'{e}tries parfaites, types de blocs, cat\'{e}gories
  d\'{e}riv\'{e}es.
\newblock {\em Ast\'{e}risque}, (181-182):61--92, 1990.

\bibitem{coleman1966idempotents}
D.~B. Coleman.
\newblock Idempotents in group rings.
\newblock {\em Proc. Amer. Math. Soc.}, 17:962, 1966.

\bibitem{linckelmann2009trivial}
Markus Linckelmann.
\newblock Trivial source bimodule rings for blocks and {$p$}-permutation
  equivalences.
\newblock {\em Trans. Amer. Math. Soc.}, 361(3):1279--1316, 2009.

\bibitem{Li18a}
Markus Linckelmann.
\newblock {\em The block theory of finite group algebras. {V}ol. {I}},
  volume~91 of {\em London Mathematical Society Student Texts}.
\newblock Cambridge University Press, Cambridge, 2018.

\bibitem{Li18b}
Markus Linckelmann.
\newblock {\em The block theory of finite group algebras. {V}ol. {II}},
  volume~92 of {\em London Mathematical Society Student Texts}.
\newblock Cambridge University Press, Cambridge, 2018.

\bibitem{nagao2014representations}
Hirosi Nagao and Yukio Tsushima.
\newblock {\em Representations of finite groups}.
\newblock Academic Press, Inc., Boston, MA, 1989.
\newblock Translated from the Japanese.

\bibitem{rickard1996splendid}
Jeremy Rickard.
\newblock Splendid equivalences: derived categories and permutation modules.
\newblock {\em Proc. London Math. Soc. (3)}, 72(2):331--358, 1996.

\bibitem{rouquier1998cyclic}
Rapha\"{e}l Rouquier.
\newblock The derived category of blocks with cyclic defect groups.
\newblock In {\em Derived equivalences for group rings}, volume 1685 of {\em
  Lecture Notes in Math.}, pages 199--220. Springer, Berlin, 1998.

\end{thebibliography}


\end{document}